\documentclass[12pt]{article}
\usepackage[cp1251]{inputenc}
\usepackage[english]{babel}
\usepackage{amssymb,amsmath}
\newtheorem{Th}{Theorem}

\newtheorem{Lem}{Lemma}
\newtheorem{Prop}{Proposition}
\newtheorem{Cor}{Corollary}

\newtheorem{Ex}{Example}
\newenvironment{Proof}
{\par\noindent{\bf Proof.}}
{\hfill$\scriptstyle\blacksquare$}

\title{Geometrical properties of the space of idempotent probability measures}
\author{A. A. Zaitov, Kh. F. Kholturaev}

\begin{document}

\maketitle
\thispagestyle{empty}

\begin{center}

\end{center}
\begin{abstract}
We establish some geometrical properties of the space of idempotent probability measures. In particular, for a compact $X$ it is established that if the space $I_{3}(X)\backslash X$ is hereditary normally, then $X$ is metrizable; some subsets allocate of the space of idempotent probability measures which are, respectively, $Z$-sets,  $\max$-$\mbox{plus}$-convex subsets, $G_\delta$-sets.\\


\end{abstract}

\tableofcontents

\section{Introduction}

The theory of idempotent measures belongs to idempotent mathematics, i. e.
the fields of the mathematics based on replacement of usual arithmetic operations
idempotent (as, for example, $x\oplus y=\max\{ x, y\} $). The notion of idempotent (Maslov) measure finds important applications in different part of mathematics, mathematical physics and economics (see the survey article ~\cite{litv2005Masldequan} and the bibliography therein). Topological and categorical properties of the functor of idempotent measures were studied in \cite{zar2010spacemapidmeasures}. Although idempotent measures are not additive and corresponding functionals are not linear, there are some parallels between topological properties of the functor of probability measures and the functor of idempotent measures (see for example \cite{zar2010spacemapidmeasures} and \cite{rad2016equiconnmonad}) which are based on existence of natural equiconnectedness structure on both functors.

However, some differences appear when the problem of the isomorphism of the functors of probability and of idempotent probability  measures was studying.

In the present paper for a compact $X$ we are established that if the space $I_{3}(X)\backslash X$ is hereditary normally, then $X$ is metrizable. Further we indicate such subsets of the space of idempotent probability measures which are, respectively, $Z$-set,  $\max$-$\mbox{plus}$-convex subset, $G_\delta$-set, $Q$-manifolds.

In the present paper under compact (\textit{pl.} compacts) we mean a compact Hausdorff space, under compactum (\textit{pl.} compacta) a metrizable compact space and under map a continuous map.

\section{On Homeomorphism of the Spaces of Probability and of Idempotent Probability measures on Compact Metrizable Space}

We will consider set $\mathbb {R} _ {\max } = \mathbb { R}\bigcup \left\{ - \infty\right\} $ with two algebraic operations:
addition $\oplus $ and multiplication $\odot$ determined as follows $u\oplus v=\max \{ u, v\} $ and $u\odot v=u+v$ where $\mathbb {R}$ is the set of real numbers.

Let $X$ be a compact, $C(X)$ be the Banach algebra of continuous functions on $X$ endowed with the usual algebraic operations and the sup-norm. For $C(X)$ operation $\oplus$ and $\odot$ we will determine as $\varphi \oplus \psi =\max \{\varphi, \psi \}$ and $\varphi\odot \psi=\varphi + \psi$, where $\varphi,\,\psi \in C(X)$. Recall that a functional $\mu :C(X)\to \mathbb{R}$ is called \cite{zar2010spacemapidmeasures} \textit{an idempotent probability measure} on $X$, if it satisfies the following properties:

(1) $\mu ({{\lambda }_{X}})=\lambda$  where ${{\lambda }_{X}}$ is a constant function on $X$ taking the value $\lambda\in
\mathbb{R}$ (normality);

(2) $\mu(\lambda \odot \varphi)=\lambda \odot \mu(\varphi)$ for all $\lambda \in \mathbb{R}$ and $\varphi \in C(X)$ (homogeneity);

(3) $\mu (\varphi \oplus \psi )=\mu (\varphi )\oplus \mu (\psi )$ for all $\varphi $, $\psi \in C(X)$ (additivity);

For a compact $X$ we denote by $I(X)$  the set of all idempotent probability measures on $X$.
Consider $ I(X)$ as a subspace of ${{\mathbb{R}}^{C(X)}}$.

For the given compacts $X$, $Y$ and a continuous mapping $f:X\to Y$  we can verify that
the natural map  $I(f):I(X)\to I(Y)$, defined by the formula $I(f)(\mu )(\psi)=\mu(\psi\circ f)$  is continuous. Moreover, the $I$ construction is a normal functor \cite{zar2010spacemapidmeasures}. Therefore, for an arbitrary idempotent probability measure $\mu \in I(X)$ one can
define the concept of \textit{support} of the measure $\mu$:
$$\text{supp}\mu =\bigcap \left\{ A\subset X:\overline{A}=A,\,\,\mu \in I(A) \right\}.$$

For a positive integer $n$ we define the following set
$${{I}_{n}}(X)=\left\{ \mu \in I(X):\,\,\left| \operatorname{supp}\mu  \right|\le n \right\}.$$

Put
$${{I}_{\omega }}(X)=\bigcup\limits_{n=1}^{\infty }{{{I}_{n}}(X)}.$$

The set ${{I}_{\omega }}(X)$ is every where dense \cite{zar2010spacemapidmeasures} in $I(X)$. An idempotent probability measure
$\mu \in {{I}_{\omega}}(X)$ is called an idempotent probability measure with finite support. Note that if $\mu$ is an idempotent probability measure with the finite support $\operatorname{supp}\mu =\left\{{{x}_{1}},{{x}_{2}},...,{{x}_{k}} \right\}$ then $\mu$ can be represented as $\mu=\lambda_{1} \odot \delta_{x_{1}} \oplus\lambda_{2} \odot \delta_{x_{2}} \oplus...\oplus \lambda_{k}\odot\delta_{x_{k}}$ uniquely, where $-\infty <{{\lambda}_{i}}\le 0$, $i=1,\dots..., k$, ${{\lambda }_{1}}\oplus {{\lambda}_{2}}\oplus ...\oplus {{\lambda }_{k}}=0$. Here, as usual, for $x\in X$ across ${{\delta }_{x}}$ we denote a functional on $C(X)$ defined by the formula ${{\delta }_{x}}(\varphi )=\varphi (x)$, $\varphi \in C(X)$, and called the Dirac measure. It is supported at the point $x$.

For a given compact $X$ the set of all probability measures, that is the set of all functionals $\mu :C(X)\to \mathbb{R}$, satisfying the conditions:

(1) $\mu ({{\lambda }_{X}})=\lambda $ for all $\lambda \in \mathbb{R}$, where ${{\lambda }_{X}}$ -- constant function;

(2) $\mu(\lambda\varphi)=\lambda\mu(\varphi)$ for all $\lambda \in \mathbb{R}$ and $\varphi \in C(X)$;

(3) $\mu (\varphi + \psi )=\mu (\varphi )+ \mu (\psi )$ for all $\varphi $, $\psi \in C(X)$,\\
denoted by $P(X)$.

For a compact $X$ the set $P(X)$ of all probability measures is endowed with the topology of pointwise convergence, i. e. we consider $P(X)$ as a supace of ${{\mathbb{R}}^{C(X)}}$.

It is well known the topological spaces $P(X)$ and $I(X)$ equipped with the topology of pointwise convergence are compacts (Hausdorff spaces).

\begin{Th}
For an arbitrary finite compact set $X$ the spaces $P(X)$ and $I(X)$ are homeomorphic.
\end{Th}

\begin{Proof} Consider the mapping $$z^{P}_{I}:P(X) \longrightarrow {I(X)},$$
given by equality
$$z^{P}_{I}\left(\sum_{i=1}^n\alpha_i\delta_{x_i}\right)=\bigoplus_{i=1}^n\left(\left(\ln \alpha_{i}-
\bigoplus_{j=1}^n \ln \alpha_{j}\right)\odot{\delta_{x_i}}\right),\qquad\sum_{i=1}^n\alpha_i\delta_{x_i}\in P(X),$$
and the mapping
$$z^{I}_{P}:I(X)\longrightarrow P(X),$$
defined by the rule
$$z^{I}_{P}\left(\bigoplus_{i=1}^n\lambda_i\odot\delta_{x_i}\right)=\sum_{i=1}^n\frac{e^{\lambda_i}}{\sum\limits_{j=1}^n e^{\lambda_j}}\cdot\delta_{x_i},\qquad\bigoplus_{i=1}^n\lambda_i\odot\delta_{x_i}\in I(X).$$

We will show that the mappings $z^{P}_{I}$ and $z^{I}_{P}$ are continuous and mutually inverse.

1) For each probability measure $\sum\limits_{i=1}^n\alpha_i\delta_{x_i}\in P(X)$ the following equalities hold $$z^{I}_{P}\left(z^{P}_{I}\left(\sum_{i=1}^n\alpha_i\delta_{x_i}\right)\right)=z^{I}_{P}\left(\bigoplus_{i=1}^n\left(\ln \alpha_i-
\bigoplus_{j=1}^n \ln \alpha_j\right)\odot\delta_{x_i}\right)=$$
$$=\sum_{i=1}^n\frac{e^{\ln \alpha_i-
\bigoplus\limits_{j=1}^n \ln \alpha_j}}{\sum\limits_{l=1}^n e^{\ln \alpha_l-\bigoplus\limits_{j=1}^n \ln \alpha_j}}\cdot\delta_{x_i}=
{\sum_{i=1}^{n}\frac{e^{\ln \alpha_{i}}:{e^{\bigoplus\limits_{j=1}^{n}\ln \alpha_{j}}}}
{\left(\sum\limits_{l=1}^{n}e^{\ln \alpha_{l}}\right):{e^{\bigoplus\limits_{j=1}^{n}\ln \alpha_{j}}}}\cdot\delta_{x_{i}}}=$$
$$={\sum_{i=1}^{n}\frac{\alpha_{i}\delta_{x_{i}}}{\sum\limits_{l=1}^{n}\alpha_{l}}}={\sum_{i=1}^{n}\alpha_{i}\delta_{x_{i}}};$$

2) For each idempotent probability measure $\bigoplus_{i=1}^{n}\lambda_{i}\odot\delta_{x_{i}}\in I(X)$ we have $$z^{P}_{I}\left(z^{I}_{P}\left(\bigoplus_{i=1}^{n}\lambda_{i}\odot\delta_{x_{i}}\right)\right)=
z^{P}_{I}\left(\sum_{i=1}^{n}\frac{e^{\lambda_{i}}}{\sum\limits_{j=1}^{n}e^{\lambda_{j}}}{\delta_{x_{i}}}\right)=$$
$$={\bigoplus_{i=1}^{n}\left(\ln \frac{e^{\lambda_{i}}}{\sum\limits_{j=1}^{n}e^{\lambda_{j}}} -
\bigoplus_{l=1}^{n}\ln \frac{e^{\lambda_{l}}}{\sum\limits_{j=1}^{n}e^{\lambda_{j}}}\right)\odot{\delta_{x_{i}}}}=$$
$$={\bigoplus_{i=1}^{n}\left(\ln e^{\lambda_{i}} - {\ln \sum_{j=1}^{n}e^{\lambda_{j}}} - \bigoplus_{l=1}^{n}
\left(\ln e^{\lambda_{l}} - \ln\sum_{j=1}^{n}e^{\lambda_{j}}\right)\right)\odot{\delta_{x_{i}}}}=$$
$$={\bigoplus_{i=1}^{n}\left(\lambda_{i}-\ln \sum_{j=1}^{n}e^{\lambda_{j}}-\bigoplus_{l=1}^{n}\lambda_{l}+
\ln \sum_{j=1}^{n}e^{\lambda_{j}}\right)\odot{\delta_{x_{i}}}}={\bigoplus_{i=1}^{n}{\lambda_{i}\odot\delta_{x_{i}}}}.$$

Consequently, the compositions $z^{P}_{I}z^{I}_{P}:I(X)\rightarrow I(X)$ and $z^{I}_{P}z^{P}_{I}:P(X)\rightarrow P(X)$ are the identical mappings.

Now we will show that the mappings $z^{P}_{I}$ and $z^{I}_{P}$ are continuous.
Since they are mutually inverse mappings between compacts, it suffices to show the continuity only of one of them.

We show that the mapping  $z^{P}_{I}:P(X)\rightarrow I(X)$  is continuous.
Let $\mu_0\in P(X)$ be a probability measure, $\{\mu_t\}_{t=1}^{\infty}\subset P(X)$ be a sequence converging to $\mu_0$ in topology of pointwise convergence (symbolically $\lim\limits_{t\rightarrow \infty} \mu_t=\mu_0$).

Since $X$ is finite, without loss of generality we can assume that $\mbox{supp}\mu_t=\mbox{supp}\mu_0=\{x_1,\ x_2,\ ...,\ x_n\}$, $n\in \mathbb{N}$, for all $t=1,\ 2,\dots$. Then we have $\lim\limits_{t\rightarrow\infty}\alpha^i_t=\alpha^i_0$, where $\mu_t=\sum\limits_{i=1}^n\alpha^i_t\delta_{x_i}$, $t=0,\ 1,\ 2,\dots$. On the other hand, the continuity of the function $\ln$ and the operation $\oplus$ implies the equality $\lim\limits_{t\rightarrow\infty}\ln \alpha^i_t=\ln \alpha^i_0$. Hence,
$$\lim\limits_{t\rightarrow\infty}\left(\ln \alpha^i_t-\bigoplus\limits_{j=1}^n \ln \alpha^j_t\right)=
\ln \alpha^i_0-\bigoplus\limits_{j=1}^n \ln \alpha^j_0.$$

Therefore, $\lim\limits_{t\rightarrow\infty}z^{P}_{I}(\mu_t)=z^{P}_{I}(\mu_0)$, i. e. the mapping $z^{P}_{I}$ is continuous.

\end{Proof}

\begin{Cor}
For an arbitrary metrizable compact $X$ the spaces $P(X)$ and $I(X)$ are homomorphic.
\end{Cor}

\section{The Functors of Probability Measures and of Idempotent Probability Measures are not Isomorphic}

The Fubini theorem is established using the $\max$-$\mbox{plus}$ variant of the Hahn-Banach theorem.

A subset $L$ of the space $C(X)$ is called \cite{zar2010spacemapidmeasures} a $\max$-$\mbox{plus}$-linear subspace in $C(X)$, if:

1) $\lambda_X\in L$ for each $\lambda\in \mathbb{R};$

2) $\lambda\odot\varphi\in L$ for each $\lambda\in \mathbb{R}$ and $\varphi\in L;$

3) $\varphi\oplus\psi\in L$ for each  $\varphi, \psi\in L.$

\begin{Lem} (\cite{zai2014interrelation}, The $\max$-$\mbox{plus}$ variant of the Hahn-Banach theorem).
Let  $L$ be a $max$-$plus$-linear subspace in $C(X)$. Let $\mu:L\longrightarrow \mathbb{R}$ be a functional satisfying the conditions of normality, homogeneity and additivity (with $C(X)$ replaced by $L$). For an arbitrary $\varphi_0\in C(X)\backslash L$ there exists an extension of the functional $\mu$ satisfying the conditions of normality, homogeneity and additivity on the minimal $max$-$plus$-linear subspace $L^{'}$ containing $L\cup\{\varphi_0\}$.
\end{Lem}

Consider the following subset in $C(X\times Y)$:

$$C_0=\left\{\bigoplus_{i=1} ^{n} \varphi_i \odot \psi_i : \varphi_i \in C(X)\ \mbox{ and } \ \psi_i \in C(Y),\  i=1, 2, ..., n,\ n\in \mathbb{N} \right\}.$$

It is obvious that $C_0$ is a $\max$-$\mbox{plus}$-linear subspace in $C(X)$.

For every pair $(\mu, \nu)\in I(X)\times I(Y)$ we put
$$(\mu\widetilde{\otimes} \nu)\left(\bigoplus_{i=1}^{n} \varphi_i \odot \psi_i\right)=\bigoplus_{i=1}^{n} \mu(\varphi_i)\odot \nu(\psi_i).$$

\begin{Prop}
$\mu\widetilde{\otimes} \nu$ is an idempotent probability measure on $C_0.$
\end{Prop}

\begin{Proof}
Each $c\in \mathbb{R}$  can be represented as  $c_{X\times Y}=a_X \odot b_Y$, where $a,\ b\in \mathbb{R}$ and $a+b=c$.
Therefore, $(\mu\widetilde{\otimes} \nu)(c_{X\times Y})=(\mu\otimes \nu)(a_X \odot b_Y)=\mu(a)\odot \nu(b)=a\odot b=c$.

Let $\lambda\in \mathbb{R}$ and $\bigoplus_{i=1}^{n} \varphi_i \odot \psi_i\in C_0$. Then \\
$(\mu\widetilde{\otimes} \nu)\left(\lambda\odot \bigoplus_{i=1}^{n} \varphi_i \odot \psi_i\right)=
(\mu\widetilde{\otimes} \nu)\left(\bigoplus_{i=1}^{n}(\lambda\odot \varphi_i)\odot \psi_i \right)=
\bigoplus_{i=1}^{n} \mu(\lambda\odot \varphi_i)\odot \nu(\psi_i)=\bigoplus_{i=1}^{n} \lambda\odot \mu(\varphi_i)\odot \nu(\psi_i)=\lambda\odot \bigoplus_{i=1}^{n} \mu(\varphi_i)\odot \nu(\psi_i)=\lambda\odot (\mu\widetilde{\otimes} \nu)\left(\bigoplus_{i=1}^{n} \varphi_i \odot \psi_i \right)$.

Finally, let $\bigoplus_{i=1}^{n} \varphi_{1\,i} \odot \psi_{1\,i}\in C_0$ and $\bigoplus_{j=1}^{m} \varphi_{2\,j} \odot \psi_{2\,j}\in C_0$.
Then\\
$(\mu\widetilde{\otimes} \nu)\left(\bigoplus_{i=1}^{n} \varphi_{1\,i} \odot \psi_{1\,i} \oplus \bigoplus_{j=1}^{m} \varphi_{2\,j} \odot \psi_{2\,j}\right)=
(\mu\widetilde{\otimes} \nu)\left(\bigoplus \varphi_{k\,l} \odot \psi_{k\,l}\right)=
\bigoplus \mu(\varphi_{k\,l})\odot \nu(\psi_{k\,l})=\bigoplus_{i=1}^{n}
\mu(\varphi_{1\,i})\odot \nu(\psi_{1\,i}) \oplus \bigoplus_{j=1}^{m} \mu(\varphi_{2\,j})\odot
\nu(\psi_{2\,j})=(\mu\widetilde{\otimes} \nu)\left(\bigoplus_{i=1}^{n}\varphi_{1\,i} \odot \psi_{1\,i}\right)
\oplus (\mu\widetilde{\otimes} \nu)\left(\bigoplus_{j=1}^{m} \varphi_{2\,j}\odot \psi_{2\,j}\right).$

\end{Proof}

Since $C_0$ is a $\max$-$\mbox{plus}$-linear subspace in $C(X\times Y)$, according to Lemma 1 there exists
an extension of the idempotent probability measure $\mu\widetilde{\otimes}\nu$ on $C(X\times Y)$.

Put
$$\mu\otimes\nu=\bigoplus\left\{\xi\in I(X\times Y):\ \xi|_{C_0}=\mu\widetilde{\otimes}\nu\right\}.$$

Thus, we have proved the following $\max$-$\mbox{plus}$ variant of the Fubini theorem.

\begin{Th}
For every pair $(\mu,\nu)\in I(X)\times I(Y)$ there exists a unique idempotent probability measure $\mu\otimes \nu\in I(X\times Y)$ such that  $(\mu\otimes \nu)(\varphi\odot\psi)=\mu(\varphi)\odot \nu(\psi)$,  $\varphi\in C(X)$, $\psi\in C(Y)$.
\end{Th}

Although the spaces $I(X)$ and $P(X)$ are homeomorphic for metrizable compacts $X$, however as the following Example 1 shows, the constructions $P$ and $I$ form different functors from each other. To construct a corresponding example in which the functors $P$ and $I$ are not isomorphic, we need some concepts from the categorical algebra (see, \cite{skornya1991obshalg}).

Let $F,\ G:\Re_1 \longrightarrow \Re_2 - $ functors from the category $\Re_1$ to the category $\Re_2$.

The natural transformation $\alpha:F\longrightarrow G$ of the functor $F$  into the functor $G$  is called [8] such a function $\alpha:ob\Re_1\longrightarrow Mor\Re_2$, such that
$\alpha_A:=\alpha(A):F(A)\longrightarrow G(A)$, ($\alpha_A$, which means the value of the function $\alpha$  in the object $A$, is a morphism of the category $\Re_2$ from the object $F(A)$  to $G(A)$ )  and for any $\varphi:A\longrightarrow B$, (from the object $A$ to the object $B$) $B\in ob\Re_1$,  in the category $\Re_2$ a commutative diagram
\[\begin{matrix}
             F(A) & \xrightarrow{\mathop{\alpha }_{A}} & G(A)  \\
             F(\varphi )\downarrow  & {} & \downarrow G(\varphi )  \\
             F(B) & \xrightarrow{\mathop{\alpha }_{B}} & G(B)  ,\\
\end{matrix}\]
i.e. $$G(\varphi)\circ\alpha_A=\alpha_B\circ F(\varphi)$$

The morphisms, $\alpha_A, A\in ob\Re_1$ are called the components of the natural transformation $\alpha$.

Thus, $\alpha=\{\alpha_A: A\in ob\Re_1\}$.
The natural transformation $\varepsilon:F\longrightarrow G$ of the functors $F,G:\Re_1\longrightarrow \Re_2$  is called [8] (natural) isomorphism if $\varepsilon_A :F(A)\longrightarrow G(A)$ -- an isomorphism of the category $\Re_2$, that is, in $\Re_2$ there exists a morphism $\eta_A: G(A)\longrightarrow F(A)$, such that $\eta_A\circ\varepsilon_A=id_{F(A)}$ and $\varepsilon_A\circ\eta_A=id_{G(A)}$  for any object $A\in ob\Re_1$.

We now give the construction of an example in which the functors $P$ and $I$ are not isomorphic.

\begin{Ex}
\end{Ex}
Consider the sets $X=\{a, b, c\}$, $Y=\{a, b\}$, $Z=\{a, c\}$, where $a$, $b$, $c$ are different points (these sets are supplied with discrete topologies). Define the following mappings:

$$f:X\longrightarrow Y,\qquad f(a)=f(c)=a,\qquad f(b)=b,$$
$$g:X\longrightarrow Z,\qquad g(a)=g(b)=a,\qquad g(c)=c.$$

We show that there is no natural transformation of the functor $ P $ to the functor $ I $ and vice versa. For this purpose, consider the objects $X$, $Y\times Z$ from the category $Comp$ and the morphism $(f,g):X\longrightarrow Y\times Z$ from the category $Comp$.

It suffices to show the map $P((f,g))$ has a property that the map $I((f,g))$ does not possess it.

The map
\[\begin{matrix}
   X & \xrightarrow{(f,g)} & Y\times Z  \\
\end{matrix}\]
under the influence of the functor $P$ goes to the map
 \[\begin{matrix}
   P(X) & \xrightarrow{P((f,g))} & P(Y\times Z)  \\
\end{matrix}\]
and
\[\begin{matrix}
   I(X) & \xrightarrow{I((f,g))} & I(Y\times Z)  \\
\end{matrix}\]
under the influence of the functor $I$.

The morphisms $P((f,g))$ and $I((f,g))$ are determined pointwise \cite{skornya1991obshalg}, i. e.
$$P((f,g))=(P(f),P(g))$$ and $$I((f,g))=(I(f),I(g)).$$

Note that, by the Fubini theorem \cite{fedfil2006gentop}, the correspondence $(P(Y)\times P(Z))\ni(\mu,\nu)\longmapsto \mu\otimes \nu \in P(Y\times Z)$ gives an embedding
$$P(Y)\times P(Z)\hookrightarrow P(Y\times Z).$$

By the $\max$-$\mbox{plus}$ variant of the Fubini theorem (Theorem 2), the correspondence
$(I(Y)\times I(Z))\ni(\mu,\nu)\longmapsto \mu\otimes \nu \in I(Y\times Z)$ gives an embedding
$$I(Y)\times I(Z)\hookrightarrow I(Y\times Z).$$
The morphisms $(P(f),P(g))$ and $(I(f),I(g))$ are determined by the formulas
$$(P(f),P(g))(\mu)(\varphi\cdot\psi)=P(f)(\mu)(\varphi)\cdot P(g)(\mu)(\psi)$$ and $$(I(f),I(g))(\nu)(\varphi\odot\psi)=I(f)(\nu)(\varphi)\odot I(g)(\nu)(\psi),$$
where $\varphi\in C(Y)$, $\psi\in C(Z)$, $\mu\in P(X)$ and $\nu\in I(X)$, respectively.

We will show the map $(P(f),P(g)):P(X)\longrightarrow P(Y)\times P(Z)$ is an embedding.

In fact, for any
 \begin{center}
    $\mu=\alpha_1 \delta_a +\alpha_2 \delta_b +\alpha_3 \delta_c$,
\end{center}
\begin{center}
    $\nu=\beta_1 \delta_a +\beta_2 \delta_b +\beta_3 \delta_c$,
\end{center}
with positive  ${\alpha_1, \alpha_2, \alpha_3, \,\beta_1, \beta_2, \beta_3,\,\, \alpha_1+\alpha_2+\alpha_3=1,\,\,\beta_1+\beta_2+\beta_3=1}$, the following equalities take place
$${P(f)(\mu)=(\alpha_1 +\alpha_3)\delta_a +\alpha_2 \delta_b},$$ $${P(g)(\mu)=(\alpha_1 +\alpha_2)\delta_a +\alpha_3 \delta_c},$$
$${P(f)(\nu)=(\beta_1 +\beta_3)\delta_a +\beta_2 \delta_b},$$ $${P(g)(\nu)=(\beta_1 +\beta_2)\delta_a +\beta_3 \delta_c}.$$
Therefore, ${(P(f),P(g))(\mu)=(P(f),P(g))(\nu)}$  if and only if
\[\left\{ \begin{matrix}
   {{\alpha }_{1}}+{{\alpha }_{3}}={{\beta }_{1}}+{{\beta }_{3}},  \\
   {{\alpha }_{2}}={{\beta }_{2}},  \\
   {{\alpha }_{1}}+{{\alpha }_{2}}={{\beta }_{1}}+{{\beta }_{2}},  \\
   {{\alpha }_{3}}={{\beta }_{3}}.  \\
\end{matrix} \right.\]
This system has a unique solution ${\alpha_1 = \beta_1, \alpha_2 = \beta_2, \alpha_3 = \beta_3}$. Hence, $\mu=\nu$.
Thus, $\left(P(f),P(g)\right)(\mu)=\left(P(f),P(g)\right)(\nu)$  if and only if $\mu=\nu$,  i. e., $\left(P(f),P(g)\right):P(X)\rightarrow P(Y)\times P(Z)$ -- is an embedding.

We show that the map $\left(I(f),I(g)\right):I(X)\rightarrow I(Y)\times I(Z)$  is not an embedding.
In fact, for idempotent probability measures
\begin{center}
    $\mu=\lambda_1 \odot \delta_a \oplus \lambda_2 \odot \delta_b \oplus \lambda_3 \odot \delta_c$,
\end{center}
\begin{center}
    $\nu=\gamma_1 \odot \delta_a \oplus\gamma_2 \odot \delta_b \oplus\gamma_3 \odot \delta_c$,
\end{center}
with $ - \infty < \lambda_1, \lambda_2, \lambda_3 , \gamma_1, \gamma_2, \gamma_3 \leq 0 $ and $\lambda_1 \oplus \lambda_2 \oplus \lambda_3 = \gamma_1 \oplus \gamma_2 \oplus \gamma_3 =0$, there are equalities
\begin{center}
    ${I(f)(\mu)= (\lambda_1\oplus\lambda_3)\odot\delta_a\oplus\lambda_2\odot\delta_b}$,
\end{center}
\begin{center}
    ${I(g)(\mu)= (\lambda_1\oplus\lambda_2)\odot\delta_a\oplus\lambda_3\odot\delta_c}$,
\end{center}
\begin{center}
    ${I(f)(\nu)= (\gamma_1\oplus\gamma_3)\odot\delta_a\oplus\gamma_2\odot\delta_b}$,
\end{center}
\begin{center}
    ${I(g)(\nu)= (\gamma_1\oplus\gamma_2)\odot\delta_a\oplus\gamma_3\odot\delta_c}$.
\end{center}
The equality $\left(I(f),I(g)\right)(\mu)=\left(I(f),I(g)\right)(\nu)$  holds if and only if
\[\left\{ \begin{matrix}
   {{\lambda }_{1}}\oplus{{\lambda }_{3}}={{\gamma }_{1}}\oplus{{\gamma }_{3}},  \\
   {{\lambda }_{2}}={{\gamma }_{2}},  \\
   {{\lambda }_{1}}\oplus{{\lambda }_{2}}={{\gamma }_{1}}\oplus{{\gamma }_{2}},  \\
   {{\lambda }_{3}}={{\gamma }_{3}}.  \\
\end{matrix} \right.\]
This system has infinitely many solutions.
For example, for every pair of $\lambda_1$ and $\gamma_1 $ with ${- \infty < \lambda_1 \leq 0}$, ${- \infty < \gamma_1 \leq 0}$ a 6-tuple $(\lambda_1, \gamma_1, 0, 0, 0, 0)$  is its solution. The equality $\left(I(f),I(g)\right)(\mu)=\left(I(f),I(g)\right)(\nu)$ is true for this 6-tuple although $\lambda_1\neq\gamma_1$. This means that the mapping $(I(f),I(g))$ is not an embedding.
Thus, there is no natural transformation of the functor $P$ to the functor $I$, since the morphism $I((f,g))=(I(f),I(g))$  is not an embedding.

In other words, there is no natural transformation $\alpha=\{\alpha_X: X\in Comp\}$ such that the diagram
\[\begin{matrix}
   I(X) & \xrightarrow{I((f,g))} & I(Y\times Z)  \\
   \mathop{\alpha }_{X}\downarrow  & {} & \downarrow \mathop{\alpha }_{Y\times Z}  \\
   P(X) & \xrightarrow{P((f,g))} & P(Y\times Z)  \\
\end{matrix}\]
would be commutative.

\section{On a Metricise Criterion of Compacts}

The set of all nonempty closed subsets of the topological space $X$ is denoted by $\exp X$. For open subsets $U_{1},...,U_{n}\subset X$     a family of the sets of the view
$$O\langle U_{1},...,U_{n}\rangle=\{F: F\in \exp X, F\subset \bigcup_{i=1}^n U_{n}, F\cap U_{1}\neq\varnothing,..., F\cap U_{n}\neq \varnothing\}$$
forms a base of a topology on $\exp X$. This topology is called \textit{Vietoris topology}, the set $\exp X$ equipped with the Vietoris topology is called a \textit{hyperspace} of the topological space $X$.

For a compact $X$ its hyperspace $\exp X$ is a compact. For the compact $X$, the natural number $n$, the functor $F$ we put
$$F_n=\{a\in F(X):\ |\mbox{supp}a|\leq n\},$$
$$F_n^{0}=F_n(X)\setminus F_{n-1}(X),$$
where
$$\text{supp}a =\bigcap \left\{ A\subset X:\overline{A}=A,\,\,a \in F(A) \right\}$$
is a support of the element $a\in F(X)$.

\begin{Prop}
For each uncountable cardinal number ${\tau}$ the space $\exp_{3}^{0}\alpha N_{\tau}$ is not normal.
\end{Prop}

\begin{Proof} Let $\tau$ be an uncountable cardinal number. Then $N_{\tau}=\{1, 2,\dots, \alpha, ...: \alpha< \tau\}$ is infinite. That is why there exist disjoint subsets $F_{1}$ and $F_{2}$ of $N_{\tau}$ such that $F_{1}$ is uncountable and $F_{2}$ is countable. Take a point $x_{0}\in F_{1}\cup F_{2}$. Choose subsets $A_{1}$ and $A_{2}$ of the space $\exp_{3}^{0}\alpha N_{\tau}$, assuming
 $$A_{1}=\left\{\{p,x,x_0\}:x\in F_{1}\setminus\{x_{0}\}\right\},\,\,\,\,\,\,A_{2}=\left\{\{p,x^{'},x_{0}\}:x^{'}\in F_{2}\setminus\{x_{0}\}\right\},$$
where $p\in \alpha N_{\tau}\setminus N_{\tau}$. Obviously, $A_{1}\cap A_{2}=\varnothing$.

Let $F=\{x_{1},x_{2},x_{3}\} \in \exp_{3}^{0}\alpha N_{\tau}\setminus A_{1}$. The set $O(\{x_1\},\{x_2\},\{x_3\})$  is an open neighbourhood of the set $F$ wich does not intersect $A_{1}$. Hence, the set $A_{1}$  is closed in  $\exp_{3}^{0}\left( {\alpha N_{\tau}} \right)$. Similarly one can check that $A_{2}$ is closed in $\exp_{3}^{0}\left( {\alpha N_{\tau}} \right)$.

For each $x\in N_{\tau}$ we put $U_{x}=O\left\langle\alpha N_{\tau}\setminus \{x_{0},x\},\{x_{0}\},\{x\}\right\rangle\cap \exp_{3}^{0} \left( {\alpha N_{\tau}} \right)$. It is easy to see that the smallest by inclusion neighbourhoods of the sets $A_{1}$ and $A_{2}$ in $\exp_{3}^{0} \left( {\alpha N_{\tau}}\right)$ are the sets $O A_{1}=\bigcup_{x\in F_{1}} U_{x}$ and $O A_{2}=\bigcup_{x\in F_{2}} U_{x}$, respectively. For the set $\{a,b,x_{0}\}$, where $a\in F_{1}$, $b\in F_{2}$, we have
$$\{a,b,x_{0}\}\in O\left(\alpha N_{\tau}\setminus \{x_{0},a\},\{x_{0}\},\{a\}\right)\subset OA_{1},$$
$$\{a,b,x_{0}\}\in O\left(\alpha N_{\tau}\setminus \{x_{0},b\},\{x_{0}\},\{b\}\right)\subset OA_{2}.$$
This means, $OA_{1}\cap OA_{2}\neq \varnothing$.

\end{Proof}

\begin{Prop}
Let ${\tau}$ is an uncountable cardinal number. Then $I_{3}^{0}\left({\alpha N_{\tau}}\right)$ is not normal space.
\end{Prop}

\begin{Proof}
For each compact $X$ the set $\exp_3^0(X)$ is closed in $I_3^0(X)$. Besides the normality is a hereditary property for the closed subsets of the space. Therefore according to Proposition 2 the space $I_{ 3}^{0}\left({\alpha N_{\tau}}\right)$ is not normal.

\end{Proof}

The following result is a metricize criterion of compacts (Hausdorff compact spaces).

\begin{Th}
Let $X$ be a compact. If the space $I_{3}(X)\backslash X$  is hereditarily normal then $X$ is metrizable.
\end{Th}

\begin{Proof} Suppose the compact $X$  is non-metrizable. If $X$ has a unique nonisolated point then $X$ is homeomorphic to $\alpha\mathbb{N}_{\tau}$ for $\tau=|X|>\omega$. Proposition 3 implies that $I_{3}^{0} (X)$ is not normal. But according to the condition $I_{3}^0(X)$ must be normal as a subset of the hereditarily normal space $I_{3}(X)\setminus X$. We get a contradiction.

Now let $a$ and $b$ be distinguished nonisolated points of the compact $X$. There are open neighbourhoods $U$ and $V$ of points $a$ and $b$, respectively, such that $\overline{U}\bigcap \overline{V}=\varnothing$.

We consider the set $Z=\overline{U}\times\exp_2\overline{V}$ and by formula $\lambda(x, y, z)=0\odot\delta_{x}\oplus  0\odot\delta_{y}\oplus 0\odot\delta_{z}$ we define the topological embedding $\lambda: Z\rightarrow I_3(X)\setminus X$.
The result of M. Katetov \cite[Corollary 1]{katetov1948complnormCartprod} (which asserts the perfectly normality of the factor $X$ under the condition of hereditarily normality of the product of $X\times Y$) implies that the factor $\exp_2\overline{V}$ of the product $Z=\overline{U}\times \exp_2\overline{V}$ is perfectly normal. Further, applying the result of V. V. Fedorchuk \cite{fed1989kteorKatetova} (which asserts the metrizability of the compact $X$ if for a normal functor $F$ of degree $\leq 2$ the space $F(X)$ is perfectly normal) we conclude that $\overline{V}$ is metrizable. Similarly, one can show that $\overline{U}$  is metrizable. Therefore each nonisolated point of the compact $X$ has a metrizable closed neighbourhood. Hence the compact $X$ is locally metrizable. Therefore it is metrizable.

\end{Proof}

\begin{Cor}
Let $X$ be a compact set and  $n\geq 3$. If the space $I_{n} (X)\backslash X$ is hereditarily normal then $X$ is metrizable.
\end{Cor}

\section{On $Z$-Sets of the Space of Idempotent Probability Measures}

A closed set $F$ of a space $X$ is said \cite{curdob1984someappl, dobr1986Z-set} to be \textit{a $Z$-set} in $X$ if all maps of compacta into $X$ can be arbitrarily closely approximated by maps into $X \setminus F$. A countable union of $Z$-sets in $X$ is called a $\sigma$-$Z$-set in $X$. It is easy to see that if the  identity map $id_X$ of $X$ has an approximation by maps into $X \setminus F$, then all maps of compacta into $X$ can be arbitrarily closely approximated by maps into $X \setminus F$.

\begin{Th}
For an arbitrary non-empty closed subset $A$ of a compact $X$, $A\neq X$, the subspace $I(A)$ is a $Z$-set in $I(X)$.
\end{Th}

\begin{Proof}
Fix $x_{0}\in X\setminus A$  and for an arbitrary $\varepsilon$, $ 0<\varepsilon\leq 1$, we define a continuous map $f_{\varepsilon}: I(X)\rightarrow I(X)$, assuming
$$f_{\varepsilon}(\mu)=(\ln (1-\varepsilon)-\ln (1-\varepsilon)\oplus \ln \varepsilon)\odot\mu\oplus (\ln \varepsilon-\ln (1-\varepsilon)\oplus \ln \varepsilon)\odot\delta_{x_{0}},\,\,\mu\in I(X).$$

The map is defined correctly. In addition, $f_{\varepsilon}(\mu)\notin I(X)$ for any $\mu \in I(X)$, because $x_{0}\in \mbox{supp}f_{\varepsilon}(\mu)$, i. e. $(X\setminus A)\cap \mbox{supp}f_{\varepsilon}(\mu)\neq\varnothing$. Therefore, the mapping $f_{\varepsilon}$ acts from $I(X)$ to $I(X)\setminus I(A)$. Consider the distance between points $x$, $y$ of the $\max$-$\mbox{plus}$ segment $[-\infty,\ 0]$ defined by the rule $\rho (x,\ y) = |\mbox{e}^y-\mbox{e}^x|$. A distance between points $\alpha_1\odot x\oplus \beta_1\odot y$ and $\alpha_2\odot x\oplus \beta_2\odot y$ of the $\max$-$\mbox{plus}$-segment $[x,\ y] = \{\alpha\odot x\oplus\beta\odot y:\ -\infty\leq \alpha\leq 0,\ -\infty\leq \beta\leq0,\ \alpha\oplus\beta=0\}$ we define by the formula
$$
\rho (\alpha_1\odot x\oplus \beta_1\odot y,\ \alpha_2\odot x\oplus \beta_2\odot y) = |\mbox{e}^{\alpha_2}-\mbox{e}^{\alpha_1}|+|\mbox{e}^{\beta_2}-\mbox{e}^{\beta_1}|.
$$
Consider the $\max$-$\mbox{plus}$-segment
$$[\mu,\ \delta_{x_0}] = \{(\ln (1-\varepsilon)-\ln (1-\varepsilon)\oplus \ln \varepsilon)\odot\mu\oplus (\ln \varepsilon-\ln (1-\varepsilon)\oplus \ln \varepsilon)\odot\delta_{x_{0}}:\ 0\leq\varepsilon \leq 1\}$$
and calculate the distance between points
$$\mu=0\odot\mu\oplus(-\infty)\odot\delta_{x_0}$$
and
$$f_\varepsilon(\mu) = (\ln (1-\varepsilon)-\ln (1-\varepsilon)\oplus \ln \varepsilon)\odot\mu\oplus (\ln \varepsilon-\ln (1-\varepsilon)\oplus \ln \varepsilon)\odot\delta_{x_{0}}$$
of the $\max$-$\mbox{plus}$-segment $[\mu,\ \delta_{x_0}]$:
$$
\rho (\mu,\ f_\varepsilon(\mu)) = |\mbox{e}^{\ln (1-\varepsilon)-\ln (1-\varepsilon)\oplus \ln \varepsilon}-\mbox{e}^0|+|\mbox{e}^{\ln \varepsilon-\ln (1-\varepsilon)\oplus \ln \varepsilon}-\mbox{e}^{-\infty}| =
$$
$$
= |\mbox{e}^{\ln (1-\varepsilon)-\ln (1-\varepsilon)\oplus \ln \varepsilon}-1|+\mbox{e}^{\ln \varepsilon-\ln (1-\varepsilon)\oplus \ln \varepsilon}.
$$

Two cases are possible.

$1^{st}$ case. $0<\varepsilon\leq\frac{1}{2}$. Then
$$
|\mbox{e}^{\ln (1-\varepsilon)-\ln (1-\varepsilon)\oplus \ln \varepsilon}-1| + \mbox{e}^{\ln \varepsilon - \ln (1-\varepsilon)\oplus \ln \varepsilon} = 0 + \mbox{e}^{\ln \frac{\varepsilon}{1-\varepsilon}}=\frac{\varepsilon}{1-\varepsilon}.
$$

$2^{nd}$ case. $\frac{1}{2}<\varepsilon\leq 1$. Then
$$
|\mbox{e}^{\ln (1-\varepsilon)-\ln (1-\varepsilon)\oplus \ln \varepsilon}-1| + \mbox{e}^{\ln \varepsilon - \ln (1-\varepsilon)\oplus \ln \varepsilon} = \mbox{e}^{\ln \frac{1-\varepsilon}{\varepsilon}} + 1 =\frac{1-\varepsilon}{\varepsilon}+1=\frac{1}{\varepsilon} .
$$

So, the approach $\varepsilon\rightarrow 0$ implies the convergence of the idempotent probability measures $f_{\varepsilon}(\mu)$ to $\mu$ by the metric $\rho$. From here follows that under the condition $\varepsilon\rightarrow 0$ the measures $f_{\varepsilon}(\mu)$ converge to $\mu$ with respect to the topology of the space $I(X)$. Indeed, for an arbitrary function $\varphi\in C(X)$ and for each $\varepsilon$, with $0<\varepsilon<\frac{1}{1+\mbox{e}^{2\|\varphi\|}}$  we have $f_{\varepsilon}(\mu)(\varphi)=0\odot\mu(\varphi)\oplus \ln \frac{\varepsilon}{1-\varepsilon}\odot\delta_{x_{0}}(\varphi) = \mu(\varphi)\oplus (\ln \frac{\varepsilon}{1-\varepsilon}+\varphi (x_{0}))=\mu(\varphi)$. Therefore, $0=|f_{\varepsilon}(\mu)(\varphi)-\mu(\varphi)| <\varepsilon$.
In other words, the identity map  $id_{I(X)}$ is approximated by the maps $f_{\varepsilon}:I(X)\rightarrow I(X)\setminus I(A)$.

\end{Proof}

\begin{Th}
For an arbitrary compact $X$ and every $n\in \mathbb{N}$, $n<\left|X\right|$, the subspace $I_{n}(X)$ is a $Z$-set in $I(X)$. Therefore, $I_{\omega}(X)$ is a $\sigma$-$Z$-set in $I(X)$.
\end{Th}

\begin{Proof} Fix $\nu \in I(X)\setminus I_{n}(X)$. Then $|\mbox{supp}\nu|\geq n+1$. For an arbitrary $\varepsilon$, $0<\varepsilon\leq 1$, we define the map $f_{\varepsilon}:I(X)\rightarrow I(X)$  by the equality
$$f_{\varepsilon}(\mu)=(\ln (1-\varepsilon)-\ln (1-\varepsilon)\oplus \ln \varepsilon)\odot\mu\oplus (\ln \varepsilon-\ln (1-\varepsilon)\oplus \ln \varepsilon)\odot\nu, \qquad \mu\in I(X).$$

The map  $f_{\varepsilon}$  is defined correctly, and it is continuous. In addition, $|\mbox{supp}f_{\varepsilon}(\mu)|\geq n+1$ because $\mbox{supp}f_{\varepsilon}(\mu)\supset \mbox{supp}\mu \cup \mbox{supp}\nu$, for every $\mu\in I(X)$. Hence, $f_{\varepsilon}(\mu)\notin I_{n}(X)$ and the map $f_{\varepsilon}$ acts from $I(X)$ to $I(X)\setminus I_n(X)$.

For each function $\varphi\in C(X)$ and every $\varepsilon$, $0<\varepsilon<\frac{1}{1+\mbox{e}^{2\|\varphi\|}}$, we have $f_{\varepsilon}(\mu)(\varphi)=0\odot\mu(\varphi)\oplus \ln \frac{\varepsilon}{1-\varepsilon}\odot\nu(\varphi) = \mu(\varphi)\oplus (\ln \frac{\varepsilon}{1-\varepsilon}+\nu(\varphi))=\mu(\varphi)$. That is why $0=|f_{\varepsilon}(\mu)(\varphi)-\mu(\varphi)| <\varepsilon$.
So, the identity map $id_{I(X)}$ is approximated by the maps $f_{\varepsilon}:I(X)\rightarrow I(X)\setminus I_n(X)$. Hence,  $I_{n}(X)$  is a $Z$-set in $I(X)$. So, $I_{\omega}(X)=\bigcup_{i=1}^{\infty}I_{n}(X)$ becomes a $\sigma$-$Z$-set.

\end{Proof}

\section{The Space of Idempotent Measures and the Hilbert Cube}

For the functor $I$, a compact $X$ and a nonempty set $A\subset X $ we set
$$S_{I}(A) = \left\{\mu\in I(X):\ \mbox{supp}\mu \cap A \ne \varnothing\right\}.$$

It is clear that $S_{I} \left({\varnothing}\right) = \varnothing $, $S_{I} \left({X} \right) = I\left( {X} \right)$. By construction the inclusion $A \subset B$ implies $S_{I} \left( {A} \right) \subset S_{I} \left( {B} \right)$.

Note that $S_{I}(A\cap B)\subset S_{I}(A)\cap S_{I}(B)$ for subsets $A$ and $B$ of a compact $X$. But the opposite is not true. Indeed, consider the sets $A=\left\{0,1,2\right\}$ and $B=\left\{1,2,3\right\}$. Then for $\mu=\lambda_{1}\odot\delta_{0}\oplus\lambda_{2}\odot\delta_{3}$, $-\infty<\lambda_{1},\ \lambda_{2}\leq 0$, $\lambda_{1}\oplus \lambda_{2}=0$, we have $\mu\in S_{I}(A)\cap S_{I}(B)$, but $\mu \notin S_{I}(A\cap B)$.

\begin{Prop}
For a compact $X$, every open subset $U\subset X$ the equality
$$I(X\backslash U) = I(X)\backslash S_{I} (U)$$
takes place.
\end{Prop}

\begin{Cor}
For every open subset $U$ of the compact $X$ the set $S_{I} \left( {U} \right)$ is open in $I\left( {X} \right)$.
\end{Cor}

\begin{Prop}
For an arbitrary non-empty set $A$ of the compact $X$, the subset  $S_{I}(A)$ is $max$-$\mbox{plus}$-convex.
\end{Prop}

\begin{Proof}
Let $\mu_{1}$, $\mu_{2}\in S_{I}(A)$. Then $A\cap \mbox{supp}\mu_{i}\neq\varnothing$, $i=1, 2$. For every pair $\alpha$,\ $\beta$ of real numbers, $-\infty<\alpha,\ \beta\leq 0$, $\alpha \oplus \beta = 0$, we have $\mbox{supp}\left(\alpha\odot\mu_{1} \oplus\beta\odot\mu_{2}\right) = \mbox{supp}\mu_{1}\cup \mbox{supp}\mu_{2}$. Hence, $\alpha\odot\mu_{1}\oplus \beta\odot\mu_{2} \in S_{I}(A)$. This proves the $\max$-$\mbox{plus}$-convexity of the set $S_{I}(A)$.

\end{Proof}

\begin{Ex}
\end{Ex}
Let $n=\left\{0,1,2,...,n-1\right\}$ be $n$-point discrete space. According to Corollary 2 for each $A\subset\left\{0,1,2,...,n-1\right\}$ the set $S_{I}(A)$ is open in $I\left(\left\{0,1,2,...,n-1\right\}\right)$. In particular, for each $i\in \left\{0,1,2,...,n-1\right\}$ the set
$$S_{I}(i)\equiv S_{I}(\{i\})=I\left(\left\{0,1,2,...,n-1\right\}\right)\setminus I\left(\left\{0,1,2,...,i-1,i+1,...,n-1\right\}\right)$$
is open. Moreover, the intersection $\bigcap\limits_{i=1}^{n-1}S_{I}\left(\left\{i\right\}\right)$ is the interior of the compact $I\left(\left\{0,1,2,...,n-1\right\}\right)$, i. e.
$$\bigcap\limits_{i=1}^{n-1}S_{I}(i)=Int\, I\left(\left\{0,1,2,...,n-1\right\}\right).$$

\begin{Prop}
For every closed subset $A$ of the compact $X$, the set $S_{I}(A)$  is a $G_{\delta}$-set in $I(X)$.
\end{Prop}

\begin{Proof}
Represent $A$ as $A=\bigcap\limits_{n=1}^{\infty}f^{-1}\left(-\frac{1}{n},\frac{1}{n}\right)$,
where $f:X\rightarrow A$ is a continuous function such that $A=f^{-1}(0)$. Then $S_{I}(A)\subset \bigcap\limits_{n=1}^{\infty}S_{I}\left(f^{-1}\left(-\frac{1}{n},\frac{1}{n}\right)\right)$. We show the reverse inclusion. Let
$\mu\in \bigcap\limits_{n=1}^{\infty}S_{I}\left(f^{-1}\left(-\frac{1}{n},\frac{1}{n}\right)\right)$. If $\mu\notin S_{I}(A)$, then $\mbox{supp}\mu\subset X\setminus A$, and therefore, $\mbox{supp}\mu\subset X\setminus f^{-1}\left(-\frac{1}{n},\frac{1}{n}\right)$ for some $n$. This implies that $\mu\notin S_{I}\left(f^{-1}\left(-\frac{1}{n},\frac{1}{n}\right)\right)$. This contradiction shows that $S_{I}(A)=\bigcap\limits_{n=1}^{\infty}S_{I}\left(f^{-1}\left(-\frac{1}{n},\frac{1}{n}\right)\right)$.

\end{Proof}

We remind \cite{cur1985bounset} the following notations: $Q=\prod\limits_{i=1}^{\infty}[-\frac{1}{i},\frac{1}{i}]$ is the \textit{Hilbert cube}, $S=\prod\limits_{i=1}^{\infty}(-\frac{1}{i},\frac{1}{i})$ and $B(Q)=Q\setminus S$ are the \textit{pseudo-interior} and the \textit{pseudo-boundary} of the Hilbert cube, respectively. A subset $B(Q)$ of the Hilbert cube $Q$ is a \textit{boundary set} if it satisfies the conditions:

1) for every $\varepsilon>0$ there is a continuous function $\eta:Q\rightarrow Q\setminus B(Q)$ such that $d(\eta, id)<\varepsilon$;

2) the complement $Q\setminus B(Q)$  is homeomorphic to $S$ (which is homeomorphic to the Hilbert space).

\textbf{Equivalent Definition.}
Let $B$ be an everywhere dense $\sigma$-$Z$-set in $Q$. It is a \textit{boundary set} in $Q$,
if $Q\setminus B \cong S$.

Let $X_{1}\subset X_{2}\subset...$ be a tower of subsets of $X$.
The system of sets $\left\{X _{n}\right\}$ is called \cite{curdob1984someappl} a \textit{strongly universal tower for compacta} if for every map $f:A\rightarrow X$ of a compactum $A$, for every closed subset $B\subset A$ such that $f|_{B}: B\rightarrow X_{m}$ is an embedding into some $X_m$, and for every $\varepsilon >0$, there exists an embedding  $h:A\rightarrow X_{n}$ for some $n\geq m$, such that $f|_{B}=h|_{B}$ and $d(f,h)<\varepsilon$. We refer to $\bigcup\limits_{i=1}^{\infty}X_{n}$  as a \textit{skeletoid for compacta}.

If in this definition for each $n\in N$ the set $X_{n}$ is a $Z$-set in $X$, then $\bigcup\limits_{i=1}^{\infty}X_{n}$ is called the $Z$ -skeletoid for compacta. Further, if in this definition $X$ is $ANR(\mathfrak{M})$-space, then $X$ is called strongly universal for arbitrary compact sets.

Let $X$ and $Y$ be topological spaces, $A\subset X,\,\,B\subset Y$. The pair $(X,A)$ \textit{is homeomorphic to the pair} $(Y, B)$ (designation $(X,A)\cong(Y,B)$) if a homeomorphism  $f:A\rightarrow B$ can be extended to the homeomorphism $\widetilde{f}:X\rightarrow Y$.

\begin{Th}
Let $X$ be an infinite compactum, and $A_{1}\subseteq A_{2}\subseteq ...$ be closed subsets of $X$, such that $\bigcup\limits_{i=1}^{\infty}A_{i}$ is everywhere dense in $X$ and $\bigcup\limits_{i=1}^{\infty}A_{i}\neq X$.
Then $\left(I(X),\bigcup\limits_{i=1}^{\infty}I(A_{i})\right)\cong \left(Q, B(Q)\right)$.
\end{Th}

\begin{Proof}
Corollary 1 implies that $\left(I(X),\bigcup\limits_{i=1}^{\infty}P(A_{i})\right)\cong \left(P(X),\bigcup\limits_{i=1}^{\infty}P(A_{i})\right)$. Theorem 3.17 \cite {fed1991prmeasintop} states $\left(P(X),\bigcup\limits_{i=1}^{\infty}P(A_{i})\right)\cong \left(Q, B(Q)\right)$.

\end{Proof}

\begin{Cor}
Let $X$ be an infinite compactum, $A_{1}\subseteq A_{2}\subseteq...$ be closed subsets of $X$ such that $\bigcup\limits_{i=1}^{\infty}A_{i}$ is everywhere dense in $X$ and $\bigcup\limits_{i=1}^{\infty}A_{i}\neq X$.
Then $I(X)\setminus \bigcup\limits_{i=1}^{\infty} I(A_{i}) \approx l_{2}$.
\end{Cor}

The following statement follows from Theorems 4 and 5.

\begin{Cor}
For an arbitrary infinite compactum $X$, the pair $\left(I(X),I_{\omega}(X) \right)$ is homeomorphic to the pair $\left (Q,B(Q) \right)$ .
\end{Cor}

Based on example 2, one can see that $\bigcap\limits_{i=1}^{\infty}S_{I}\left(Q\setminus W_{i}^{\pm}\right)\cong S$. Since $S\subset Q\setminus W_{i}^{\pm}$ for each $i=1,2,...$, then $$S_{I}(S)\subset \bigcap\limits_{i=1}^{\infty}S_{I}\left(Q\setminus W_{i}^{\pm}\right).$$
Proposition 6 implies that
$$I(BdQ)=I\left(\bigcup\limits_{i=1}^{\infty}W_{i}^{\pm}\right)=I(Q\setminus S)=I(Q)\setminus S_{I}(S).$$
therefore
$$I(BdQ)\supset \bigcup\limits_{i=1}^{\infty}\left(I(Q)\setminus S_{I}\left(Q\setminus W_{i}^{\pm}\right)\right)=\bigcup\limits_{i=1}^{\infty}I\left(W_{i}^{\pm}\right).$$

Consequently, $I(BdQ)$ contains $Z$-skeletoid.
Since $I(BdQ)=\bigcup\limits_{i=1}^{\infty}S_{I}\left(W_{i}^{\pm}\right)$, then $S_{I}(BdQ)$ also contains
$Z$-skeletoid. For every $Z$-skeletoid  $\bigcup\limits_{n=1}^{\infty}A_{n}\subset Q,\,\,A_{n}\cap S\neq \varnothing,\,\, n=1,2,...$ the set $\bigcup\limits_{n=1}^{\infty}I(A_{n})$  is a $Z$-skeletoid in $S_{I}(S)$. So, we have established

\begin{Th}
The spaces $S_{I}(BdQ),\,\,S_{I}(S)$ and $I(BdQ)$ contain the skeletoid for compacta.
\end{Th}

\end{document}